\documentclass[a4paper,11pt,reqno]{amsart}
\usepackage[latin1]{inputenc}
\usepackage[english]{babel}
\usepackage{amssymb,amsfonts, amsmath, color}
\usepackage[margin=2.7cm]{geometry}
\usepackage{xcolor}
\usepackage{graphicx, color, enumerate}
\usepackage[latin1]{inputenc}
\usepackage[active]{srcltx}
\usepackage{tikz}
\usepackage{pgf}
\usepackage{etex}
\usepackage{verbatim}
\usepackage{tikz-3dplot}
\usepackage{pgfkeys}
\usepackage{amsmath}
\usepackage{amsfonts}
\usepackage{amssymb}
\usepackage{amsthm}
\usepackage{algpseudocode}
\usepackage{float}
\usepackage{xcolor}
\usepackage{mathrsfs}
\usepackage{import}
\usepackage{geometry}
\usepackage{fancyhdr}
\usepackage{fp}
\usepackage[colorlinks, citecolor=blue, linkcolor=red]{hyperref}

\newtheorem{theorem}{Theorem}[section]
\newtheorem{proposition}[theorem]{Proposition}
\newtheorem{corollary}[theorem]{Corollary}
\newtheorem{lemma}[theorem]{Lemma}
\newtheorem{remark}[theorem]{Remark}
\newtheorem{definition}[theorem]{Definition}

\newcommand{\bcl}{\begin{center}}
\newcommand{\ecl}{\end{center}}
\newcommand{\brl}{\begin{right}}
\newcommand{\erl}{\end{right}}
\newcommand{\ben}{\begin{enumerate}}
\newcommand{\een}{\end{enumerate}}
\newcommand{\overliner}{\begin{array}}
\newcommand{\earr}{\end{array}}
\newcommand{\btab}{\begin{tabular}}
\newcommand{\etab}{\end{tabular}}
\newcommand{\bdoc}{\begin{document}}
\newcommand{\edoc}{\end{document}}
\newcommand{\beqy}{\begin{eqnarray}}
\newcommand{\eeqy}{\end{eqnarray}}

\newcommand{\beqi}{\begin{eqnarray*}}
\newcommand{\eeqi}{\end{eqnarray*}}
\newcommand{\bitem}{\begin{itemize}}

\newcommand{\eitem}{\end{itemize}}
\newcommand{\nln}{\newline}
\newcommand{\newt}{\newtheorem}

\newcommand{\pa}{\partial}
\newcommand{\re}{{I\!\!R}}
\newcommand{\Rn}{\R^N}
\newcommand{\xr}{x\in\R }
\newcommand{\x}{\times}
\newcommand{\dyle}{\displaystyle}
\newcommand{\ene}{{I\!\!N}}
\newcommand{\irn}{\int\limits_{\R^N}}
\newcommand{\io}{\int\limits_{\O}}
\newcommand{\meas}{{\rm meas\,}}

\newcommand{\sign}{{\rm sign}}
\newcommand{\map}{\longrightarrow }
\newcommand{\imp}{\Longrightarrow }
\renewcommand{\div}{\nabla\cdot }
\newcommand{\sen}{{\rm sen\,}}
\newcommand{\tg}{{\rm tg\,}}
\newcommand{\arcsen}{{\rm arcsen\,}}
\newcommand{\arctg}{{\rm arctg\,}}
\newcommand{\supp}{{\textsl supp\ }}
\newcommand{\ity}{\int_{-\iy}^{+\iy}}
\newcommand{\limit}{\lim\limits}
\newcommand{\limi}{\limit_{n\to\infty}}
\newcommand{\sumi}{\sum\limits_{n=1}^{\infty}}
\newcommand{\ulu}{\underline u}
\newcommand{\ulw}{\underline w}
\newcommand{\ulz}{\underline z}
\newcommand{\ulv}{\underline v}
\newcommand{\uls}{\underline s}
\newcommand{\olu}{\overline u}
\newcommand{\olv}{\overline v}
\newcommand{\ols}{\overline s}
\newcommand{\ob}{\overline\b}
\newcommand{\ovar}{\overline\var}
\newcommand{\wv}{\widetilde v}
\newcommand{\wu}{\widetilde u}
\newcommand{\ws}{\widetilde s}
\renewcommand{\a }{\alpha }
\renewcommand{\b }{\beta }
\newcommand{\g }{\gamma}
\newcommand{\G }{\Gamma }
\renewcommand{\d }{\delta }

\newcommand{\D }{\Delta }
\newcommand{\e }{\varepsilon }
\newcommand{\z }{\zeta }
\renewcommand{\l }{\lambda }
\renewcommand{\L }{\Lambda }
\newcommand{\m }{\mu }
\newcommand{\n }{\nabla }
\newcommand{\s }{\sigma }
\newcommand{\Sig }{\Sigma }
\renewcommand{\t }{\tau }
\newcommand{\var }{\varphi }
\renewcommand{\o }{\omega }
\renewcommand{\O }{\Omega }
\newcommand{\R}{{\mathbb{R}}}
\newcommand{\bC}{{\bf C}}
\newcommand{\bZ}{{\bf Z}}
\newcommand{\bN}{{\bf N}}
\newcommand{\bQ}{{\bf Q}}
\newcommand{\bK}{{\bf K}}
\newcommand{\bI}{{\bf I}}
\newcommand{\bv}{{\bf v}}
\newcommand{\bV}{{\bf V}}
\DeclareMathOperator{\suppo}{supp} \DeclareMathOperator{\di}{div}

\newenvironment{Proof}{\Rmovelastskip\vskip12pt
plus 1pt \noindent\em\rm}{\hfill {\qed \hskip .2cm}}

\begin{document}

\title[Uniqueness for an elliptic equation with drift on manifolds]{Uniqueness in weighted Lebesgue spaces for an elliptic equation with drift on manifolds}

\author{Giulia Meglioli}

\address{\hbox{\parbox{5.7in}{\medskip \noindent{Giulia Meglioli, \\Dipartimento di Matematica, \\Politecnico di Milano, \\Piazza Leonardo da Vinci 32, 20133, Milano, Italy \\ [3pt] \emph{E-mail address: }{\tt giulia.meglioli@polimi.it}}}}}

\author{Alberto Roncoroni}

\address{\hbox{\parbox{5.7in}{\medskip \noindent{Alberto Roncoroni, \\Dipartimento di Matematica, \\Politecnico di Milano, \\Piazza Leonardo da Vinci 32, 20133, Milano, Italy \\ [3pt] \emph{E-mail address: }{\tt alberto.roncoroni@polimi.it}}}}}



\maketitle

\begin{abstract}
We investigate the uniqueness, in suitable weighted Lebesgue spaces, of solutions to a class of elliptic equations with a drift posed on a complete, noncompact, Riemannian manifold $M$ of infinite volume and dimension $N\ge2$. Furthermore, in the special case of a model manifold with polynomial volume growth, we show that the conditions on the drift term are sharp.

\end{abstract}

\
\section{Introduction}\setcounter{equation}{0}
\

We are concerned with uniqueness of solutions, in suitable weighted Lebesgue spaces, to the following linear elliptic equation
\begin{equation}\label{problema}
\Delta u + g( b,\,\nabla u ) -c\,u=0 \quad  \text{in}\,\,\,M,
\end{equation}
where $M$ is an $N$-dimensional, $N\geq 2$, complete, noncompact, Riemannian manifold of infinite volume, endowed with a metric tensor $g$. Here, $\Delta$ and $\nabla$ denote, respectively, the Laplce-Beltrami operator and the gradient on $M$ with respect to $g$; $b:M\to M$ is a given vector field and $c$ is a nonnegative function defined in $M$. For all vector fields $X,Y$ belonging to the tangent space of $M$, we set
$$
\left\langle X, Y\right\rangle:=g(X,Y)\,.
$$
\noindent Moreover, for any $x_0\in M$ and $R > 0$, let $B_R(x_0) := \{x\in M:\,\, \operatorname{dist}(x, x_0) < R\}$ denote the geodesic ball of radius $R$ and centred at $x_0$, where $\operatorname{dist}(x, x_0)$ is the geodesic distance between $x$ and $x_0$. Furthermore, let $V(x_0, R)$ denote the Riemannian volume of $B_R(x_0)$. In what follows, we set
\begin{equation}\label{r}
r(x)\equiv r:=\operatorname{dist}(x,x_0).
\end{equation}
Let us define the following subsets of $M$
\begin{equation}\label{d}
\begin{aligned}
&D_+:=\{x\in M\,\,:\, \langle b(x),x\rangle> 0\};\\
&D_-:=\{x\in M\,\,:\, \langle b(x),x\rangle\le 0\}.
\end{aligned}
\end{equation}
Then, concerning the vector field $b$, we will assume either that:
\begin{equation}\tag{{\it $H_0$}}\label{h0}
\begin{aligned}
&\textrm{(i)} \; b:M\to M,\,\,\,b\in C^1(M)\,;\\
&\textrm{(ii)} \;\text{there exist}\,\,\,K_1>0\,\,\text{such that}\,\,\\
&\quad\quad \quad |b(x)|\le K_1\quad\text{for all}\,\,x\in D_+;\\
&\textrm{(iii)}\;[\operatorname{div} b(x)]_{-}\le K_1\quad\text{for all}\,\,x\in M;
\end{aligned}
\end{equation}
or that, for some $0<\theta<1$ that will be specified later,
\begin{equation}\tag{{\it $H_1$}}\label{h4}
\begin{aligned}
&\textrm{(i)} \; b:M\to M,\,\,\,b\in C^1(M)\,;\\
&\textrm{(ii)} \;\text{there exist}\,\,\sigma\le 1-\theta\,\,\text{and}\,\,K_3>0\,\,\text{such that}\,\,\\
&\quad\quad \quad \langle b(x), \nabla r\rangle\le K_3(1+r)^{\sigma}\quad\text{for all}\,\,x\in D_+;\\
&\textrm{(iii)}\;[\operatorname{div} b(x)]_{-}\le K_3(1+r)^{\sigma-1}\quad\text{for all}\,\,x\in M;
\end{aligned}
\end{equation}
or that
\begin{equation}\tag{{\it $H_2$}}\label{h1}
\begin{aligned}
&\textrm{(i)} \; b:M\to M,\,\,\,b\in C^1(M)\,;\\
&\textrm{(ii)} \;\text{there exist}\,\,\sigma\le 1\,\,\text{and}\,\,K_2>0\,\,\text{such that}\,\,\\
&\quad\quad \quad \langle b(x), \nabla r\rangle\le K_2(1+r)^{\sigma}\quad\text{for all}\,\,x\in D_+;\\
&\textrm{(iii)}\;[\operatorname{div} b(x)]_{-}\le K_2(1+r)^{\sigma-1}\quad\text{for all}\,\,x\in M.
\end{aligned}
\end{equation}

\noindent Here, for a given function $v$, the negative part $[\cdot]_{-}$ is defined as
$$
[v]_{-}:=\max\left\{0;\,-v\right\}.
$$
Moreover, the coefficient $c$ is such that
\begin{equation}\label{h2}\tag{{\it $H_3$}}
c\in C(M),\;\,\, c(x)\geq c_0 \quad \textrm{for all}\;\; x\in M\,,
\end{equation}
for some $c_0>0$.
\medskip

\noindent 
We now summarize our main results and we refer to Section \ref{sec1} for the precise statements.

\subsection{Outline of our results} 
The main results of this paper will be given in detail in the forthcoming Theorems \ref{teo1}, \ref{teo3} and \ref{teo2}. We give here a sketchy outline of these results, describing motivations and techniques of proofs.
\medskip
For any $\phi\in C(M)$, $\phi>0$, $p> 1$, set
$$
L^p_\phi(M):= \left\{u: M\to \R\,\;\textrm{measurable}\;\,|\, \int_{M}|u|^p\phi(x) dx <\infty\right\}\,.
$$

\medskip

\noindent We shall prove the following 

\begin{itemize}

\item Let $M$ be a complete, noncompact, Riemannian manifold of infinite volume such that, for some $x_0\in M$,
\begin{equation}\label{vol1}
V(x_0,r)\le e^{\alpha r},\quad \text{for some}\,\,\,\alpha>0.
\end{equation}
If $b:M\to M$ satisfies assumption \eqref{h0}, then the solution to equation \eqref{problema} is unique in the class $L^p_{\psi}(M)$ with $p>1$ and
\begin{equation}\label{eq12}
\psi(x)=\psi(r):=e^{-\beta\,r}\;\, (r\in (0,+\infty))\,,
\end{equation}
for properly chosen $\beta>\alpha$ (see Theorem \ref{teo1}). Observe that, in general, our uniqueness class includes {\it unbounded} solutions. Thus, in particular, we get uniqueness of bounded solutions. More precisely, bounded solutions always belong to this class of uniqueness $L^p_{\psi}(M)$. Furthermore, also unbounded solutions can be considered if such solutions satisfy
$$
u(x)\le Ce^{\gamma r}\quad\text{for}\,\,\,r>1,\,\,\gamma>0\,,
$$
provided that $\beta>\alpha+\gamma p$.

\item Let $M$ be a complete, noncompact, Riemannian manifold of infinite volume such that, for some $x_0\in M$,
\begin{equation}\label{vol3}
V(x_0,r)\le e^{\alpha r^{\theta}},\quad \text{for some}\,\,\,\alpha>0,\quad \text{for some}\,\,\,0<\theta<1.
\end{equation}
If $b:M\to M$ satisfies assumption \eqref{h4} with $\theta$ as in \eqref{vol3}, then the solution to equation \eqref{problema} is unique in the class $L^p_{\eta}(M)$ with $p>1$ and
\begin{equation}\label{eq2222}
\eta(r):=e^{-\beta\,r^{\theta}} \quad \text{for some}\,\,\,\beta>0,
\end{equation}
for properly chosen $\beta>\alpha$ (see Theorem \ref{teo3}). Observe that, in particular, bounded solutions always belong to this class of uniqueness $L^p_{\eta}(M)$. Moreover, also unbounded solutions can be considered if such solutions satisfy
$$
u(x)\le Ce^{\gamma r^{\theta}}\quad\text{for}\,\,\,r>1,\,\,\gamma>0\,,
$$
provided that $\beta>\alpha+\gamma p$.

\item Let $M$ be a complete, noncompact, Riemannian manifold of infinite volume such that, for some $x_0\in M$,
\begin{equation}\label{vol2}
V(x_0,r)\le r^{\alpha},\quad \text{for some}\,\,\,\alpha>0.
\end{equation}
If $b:M\to M$ satisfies assumption \eqref{h1}, then the solution to equation \eqref{problema} is unique in the class $L^p_{\xi}(M)$ with $p>1$ and
\begin{equation}\label{eq13}
\xi(x)=\xi(r):=(1+r)^{-\tau}\;\, (r\in (0,+\infty))\,,
\end{equation}
for properly chosen $\tau>\alpha+N-1$ (see Theorem \ref{teo2}). Similarly to the previous case, we observe that, our uniqueness class includes {\it unbounded} solutions. Thus, in particular, we get uniqueness of bounded solutions. More precisely, bounded solutions always belong to this class of uniqueness $L^p_{\xi}(M)$. Moreover, also unbounded solutions can be considered if such solutions satisfy
$$
u(x)\le Cr^{\gamma}\quad\text{for}\,\,\,r>1,\,\,\gamma>0\,,
$$
provided that $\beta>\alpha+\gamma p+N-1$.

\end{itemize}

\noindent Observe that hypotheses \eqref{h4} and \eqref{h1} allow $b$ to be {\it unbounded}. In addition, such assumptions require a growth condition on $b$, when $\langle b(x), x\rangle> 0$. However, when $\langle b(x), x\rangle\leq 0$ no further conditions on $b$ are imposed. 

\begin{itemize}

\item Furthermore, we show sharpness of the result of Theorem \ref{teo2} in the special case of a model manifold $M$ such that \eqref{metric} holds. More precisely, we show that, if the drift term $b$ violetes, in an appropriate sense, assumption \eqref{h1} (that is, if it satisfies inequality \eqref{h3} below), then infinitely many bounded solutions to problem \eqref{problema} exist (see Proposition \ref{prop4} and Corollary \ref{cor} below).
\end{itemize}


\noindent 
We shall now briefly recall some known results related to equation \eqref{problema}.

\subsection{Summary of known results}\label{pippo}
Let us firstly consider the Euclidean setting. In \cite{AB} uniqueness of solutions to the Cauchy problem
 \begin{equation}\label{eq18}
\begin{cases}
\partial_t u = L\, u\quad &\textrm{in}\,\,\R^N\times(0,T)\\
u \, = 0 \quad&\textrm{in}\,\, \R^N\times \{0\} \,,
\end{cases}
\end{equation}
is investigated. Here $T>0$ and 
$$
L u:= \sum_{i,j=1}^N\frac{\pa^2 \big[a_{ij}(x,t) u\big] }{\pa x_i\pa x_j} - \sum_{i=1}^N\frac{\pa \big[b_i(x,t) u\big]}{\pa x_i}+ c(x,t)u \,,
$$
the coefficients which appear in $L$, together with all their derivatives, are locally bounded functions and the matrix $A\equiv (a_{ij})$ is assumed to be positive semidefinite in $\R^N\times(0,T)$. Moreover the authors assume that $$
|a_{ij}(x,t)|\leq K_1 ( 1+|x|^2)^{\frac{2-\l}{2}},\; |b_i(x,t)|\leq K_2(|x|^2+1)^{\frac{1}{2}},\;|c(x,t)|\leq K_3(|x|^2+1)^{\frac{\l}2},
$$
for almost every $(x,t)\in \R^N\times(0,T)$, for some constants $\l\geq 0$ and $K_i>0$, for $i=1,2,3$.  Given $\phi\in C(\R^N\times(0,T))$, $\phi>0$, set
$$
L^1_\phi(\R^N\times(0,T)):= \left\{u: \R^N\times(0,T)\to \R\,\;\textrm{measurable}\;\,|\, \int_0^T\int_{\R^N}|u|\phi(x,t) dx dt<\infty\right\}\,.
$$
In \cite[Theorem 1]{AB} it is shown that if $u$ is a solution to problem \eqref{eq18} and $u\in L^1_{\phi}(\R^N\times(0,T))$, with
\begin{equation}\label{eq19}
\phi(x)=(|x|^2+1)^{-\a_0} \quad (x\in \R^N)\;\; \textrm{if}\;\;
\l=0,
\end{equation}
or
\begin{equation}\label{eq15}
\phi(x)=e^{-\a_0(|x|^2+1)^{\frac{2-\a}2}}\quad (x\in \R^N)\;\;
\textrm{if}\;\; \l>0\,,
\end{equation}
for some $\a_0>0$, then
$$
u\equiv 0\quad \textrm{in}\;\; \R^N\times(0,T)\,.
$$
The result in \cite{AB} has attracted a lot of interests and analogue results has been proved in \cite{EKP,IKO,KPT,PoT} (we refer to these papers for further discussions and details) and generalized to the fractional heat equation in \cite{MP, PV1,PV2}. We also mention that similar results has been considered in the case of degenerate elliptic and parabolic problems in $\mathbb{R}^N$ and also in bounded domains of $\mathbb{R}^N$ (see e.g. \cite{MontPun,PPT2,PoT,P,PTes} and the references therein for previous results).
\medskip

Let us now consider the Riemannian setting. The study of uniqueness of solutions to the heat equation
\begin{equation}\label{heat}
\begin{cases}
\partial_t u=\Delta u\quad &\text{in} \,\,M\times(0,T)\\
u=u_0\quad &\text{in}\,\, M\times\{0\},
\end{cases}
\end{equation}
where $M$ is a complete Riemannian manifolds, has been largely investigated in the literature. Indeed, it is well-known that uniqueness of solutions to problem \eqref{heat} is equivalent to the stochastic completeness of the Riemannian manifold $M$ (see e.g. \cite[Theorem 6.2]{Grig}). We refer the interested reader to \cite{Grig} for a complete and exhaustive picture about these kind of results.  More precisely, in \cite[Theorem 9.2]{Grig} it is shown that, if $u$ is a classical solution to \eqref{heat} with $u_0\equiv0$, then 
$$
u\equiv 0\quad \textrm{in}\;\; M\times(0,T)\,,
$$
provided that $u\in L^2_\phi(M\times(0,T))$, where 
$$
L^2_\phi(M\times(0,T)):= \left\{u: M\times(0,T)\to \R\,\;\textrm{measurable}\;\,|\, \int_0^T\int_{M}u(x,t)^2\phi(x) d\mu(x)\, dt <\infty\right\}\,,
$$
with 
$$
\phi(x)=e^{-f(r(x))}\, ,
$$
where $r$ is given in \eqref{r} and $f$ is a positive increasing continuous function defined in $(0,\infty)$ such that
$$
\int_{R_0}^{\infty}\frac{r}{f(r)}\, dr=\infty\, ,
$$ 
for some $R_0>0$. The previous result in \cite{Grig} has been generalized in \cite{P3} in the following sense: if $u$ is a classical solution to \eqref{heat} with $u_0\equiv 0$, then 
$$
u\equiv 0\quad \textrm{in}\;\; M\times(0,T)\,,
$$
provided $u\in L^p_\phi(M\times(0,T))$, for $1<p\leq 2$, where 
$$
L^p_\phi(M\times(0,T)):= \left\{u: M\times(0,T)\to \R\,\;\textrm{measurable}\;\,|\, \int_0^T\int_{M}\vert u(x,t)\vert^p\phi(x) d\mu(x)\, dt <\infty\right\}\,,
$$
with $\phi$ as before.

\medskip

Concerning the elliptic case, in \cite{Grig2} (see also \cite[Section 13.2]{Grig}) the author proves that the only bounded solution to equation \eqref{problema}, with $b\equiv 0$ and $c(x)\ge0$, is the trivial one provided that
$$
V(x_0,r)\le kr^2e^{kC^2(r/2)}, \quad \text{for all}\,\,\,r>0\,\,\, \text{large enough}
$$
where $k>0$ and
$$
C(r):=\int_{0}^{r}\sqrt{\inf_{x\in\partial B_s(x_0)}c(x)}\,ds\,.
$$
Similar results have also been obtained for more general linear, elliptic, degenerate equations such as
\begin{equation}\label{02}
\it{L}u-cu=f \quad\quad\text{in}\,\,\,\Omega\,;
\end{equation}
where
$$
\it{L} u:= \sum_{i,j=1}^Na_{ij}\frac{\partial^2  u }{\pa x_i\pa x_j} +\sum_{i=1}^Nb_i\frac{\pa u}{\pa x_i}\,,
$$
with possibly unbounded coefficients $a_{ij}$ , $b_i$, $c$ and function $f$.
Well-posedness of problem \eqref{02} has been intensively investigated both in the case of $\Omega=\R^N$ and $\Omega$ being a bounded domain where the problem is completed with suitable boundary conditions and the coefficients can be degenerate or singular at the boundary of the domain (see e.g. \cite{P2,PTes} and the references therein). Similar uniqueness results have been studied also in the framework of nonlocal diffusion, we refer to \cite{MP,PV1,PV2}.

In particular, in \cite[Section 1.2]{P2} the author consider the problem \eqref{02} in a bounded and regular domain $\Omega$ of $\mathbb{R}^N$ where the coefficients $a_{ij}$ and $b_i$ satisfy suitable assumptions (see \cite[$(H_1)-(H_3)-(H_4)$]{P2}), $c$ satisfies \eqref{h2} and $f\in C(\Omega)$. Under these hypothesis in \cite[Theorem 1.5]{P2} it is shown that there exist at most one solution $u\in L^p_{\phi}(\Omega)$ of problem \eqref{02}, for some suitable $p>1$. In this case, 
$$
\phi(x)=[d(x)]^{\alpha}\, 
$$
where $\alpha>0$ and $d(x)=\mathrm{dist}(x,\partial\Omega)$. Of course, if $f\equiv 0$ then one obtains a uniqueness result in the weighted Lebesgue space in the same spirit of the results that we are going to present.

\bigskip

\noindent \textbf{Organisation of the paper.} The paper is organized as follows. In Section \ref{sec1} we state our main results and we recall some basic facts from Riemannian geometry that will be useful in the sequel.  Section \ref{p1}, Section \ref{p4} and Section \ref{p2} are devoted to the proof of Theorem \ref{teo1}, Theorem \ref{teo3} and Theorem \ref{teo2}, respectively. Finally, in Section \ref{p3}, we show Proposition \ref{prop4} and Corollary \ref{cor}.

\

\section{Mathematical framework and results} \label{sec1}\setcounter{equation}{0}

\

Throughout the paper we deal with classical solutions to equation \eqref{problema}, for this reason we recall their definition.

\begin{definition}\label{defsole}
We say that a function $u$ is a classical solution to equation \eqref{problema} if
\begin{itemize}
\item[(i)]  $u\in  C^{2}(M)$;
\item[(ii)]  $\Delta u - \left \langle b(x),\,\nabla u \right \rangle + c(x)u \,=\,0$\;\;  for all\;\; $x\in M$\,.
\end{itemize}
Furthermore, we say that $u$ is a supersolution (subsolution) to equation \eqref{problema}, if in $(ii)$ instead of $``="$ we have $``\geq"\; (``\leq")$\,.
\end{definition}

\noindent We now state the main result where we consider those Riemannian manifolds $M$ such that \eqref{vol1} is satisfied and we assume the vector field $b$ to be bounded, see \eqref{h0}.

\begin{theorem}\label{teo1}
Let $M$ be a complete, noncompact, Riemannian manifold such that \eqref{vol1} is satisfied. Let assumptions \eqref{h0} and \eqref{h2} be satisfied. Let $u$ be a classical solution to equation \eqref{problema}. Let $\psi$ be as in \eqref{eq12}. Let $p>1$, with $p c_0$ is large enough. If $u\in L^p_{\psi}(M)$, then
$$
u\equiv 0 \quad \textrm{in} \;\; M\,.
$$
\end{theorem}

\noindent The next two results show that, if one restricts the assumptions on the Riemannian manifold, thus \eqref{vol1} does not follows anymore, then $b$ can also be unbounded in an appropriate sense that, somehow, is related to the volume growth condition made on $M$.

\smallskip

\begin{theorem}\label{teo3}
Let $M$ be a complete, noncompact, Riemannian manifold such that \eqref{vol3} is satisfied. Let assumptions \eqref{h4} and \eqref{h2} be satisfied. Let $u$ be a classical solution to equation \eqref{problema}. Let $\eta$ be as in \eqref{eq2222}. Let $p>1$, with $p c_0$ is large enough. If $u\in L^p_{\eta}(M)$, then
$$
u\equiv 0 \quad \textrm{in} \;\; M\,.
$$
\end{theorem}

\smallskip

\begin{theorem}\label{teo2}
Let $M$ be a complete, noncompact, Riemannian manifold such that \eqref{vol2} is satisfied. Let assumptions \eqref{h1} and \eqref{h2} be satisfied. Let $u$ be a classical solution to equation \eqref{problema}. 
Let $\xi$ be defined as in \eqref{eq13}. Let $p>1$, with $p c_0$ is large enough. If $u\in L^p_{\xi}(\R^N)$, then
$$
u\equiv 0 \quad \textrm{in} \;\; M\,.
$$
\end{theorem}

\medskip

We highlight that our theorems are new even in the case $b\equiv 0$. Indeed, as already mentioned in Section \ref{pippo}, the uniqueness of solutions of equation \eqref{problema}, with $b\equiv 0$, was known only for bounded solutions (see \cite{Grig2}).

\medskip

\noindent In the following we collect some observations regarding the hypothesis of the theorems: 

\begin{itemize}
\item[$i)$] The hypothesis $p \,c_0 $ large enough made in Theorem \ref{teo1}, Theorem \ref{teo3} and Theorem \ref{teo2} will be specified in their proofs. 
\item[$ii)$] Observe that assumption \eqref{vol3} in Theorem \ref{teo3} is stronger than assumption \eqref{vol1} in Theorem \ref{teo1}. Although, this hypothesis enables us to relax the assumption on the vector field $b$, thus $b$ in Theorem \ref{teo3} might be unbounded, see \eqref{h4}.
\item[$iii)$] We mention that thanks to the Laplacian comparison Theorems (see e.g.  \cite[Section 2]{GW} or \cite[Section 15]{Grig}) we have that assumptions \eqref{vol1}, \eqref{vol3} and \eqref{vol2} are satisfied if the Ricci curvature of the Riemannian manifold $M$ is bounded from below. In particular, \eqref{vol1} and \eqref{vol3} are guaranteed if 
$$
\mathrm{Ric}\ge -k\quad \text{for some}\,\,\, k>0\,.
$$
Similarly, \eqref{vol2} is true if
$$
\mathrm{Ric}\ge 0\,.
$$
\end{itemize}



\bigskip

\noindent Finally, in the special case of a model manifold $\mathbb{M}^N_{\varphi}$ (see Definition \ref{def_modello}) satisfying \eqref{vol2} with 
\begin{equation}\label{metric}
\varphi(r)=r^{\lambda},\quad \text{for any}\,\,\,r>1\,,\quad \text{and for some}\,\,\,\lambda<\infty \,,
\end{equation}
we show that the hypothesis \eqref{h1} on $b$, is sharp: this is the content of the next proposition.  In particular, we observe that due to \eqref{volume_balle}, assumption \eqref{metric} is compatible with \eqref{vol2}, provided $\lambda<\frac{\alpha-1}{N-1}$.
\smallskip

\medskip
\begin{proposition}\label{prop4}
Let $\mathbb{M}^N_{\varphi}$ be a model manifold such that \eqref{metric} and \eqref{vol2} are satisfied. Suppose that assumption \eqref{h2} holds and let $b\in C^1(\mathbb{M}^N_{\varphi})$ be such that
\begin{equation}\label{h3}
\begin{aligned}
\langle b(x), x \rangle \geq 0 \quad &\text{ for all }\,\, x\in \mathbb{M}^N_{\varphi}\,,\\
\langle b(x),\nabla r \rangle \ge\, K_2\, r^{\sigma}\quad&\text{ for all }\,\,x\in \mathbb{M}^N_{\varphi}\setminus B_{R_0},
\end{aligned}
\end{equation} 
for some $$R_0>1, \quad \sigma>1,\quad \text{and}\quad K>1.$$
Then equation \eqref{problema} admits infinitely many solutions in the sense of Definition \ref{defsole}.
\end{proposition}

\smallskip

\noindent In the following corollary we summarize the picture of uniqueness and non-uniqueness results in the special case of model manifolds.

\begin{corollary}\label{cor}
Let $\mathbb{M}^N_{\varphi}$ be a model manifold such that \eqref{metric} and \eqref{vol2} are satisfied. Suppose that assumption \eqref{h2} holds and let $b\in C^1(\mathbb{M}^N_{\varphi})$.
\begin{itemize}
\item[(i)] Assume that $b$ satisfies assumption \eqref{h1}. Let $u$ be a classical solution to equation \eqref{problema}. Let $p>1$ and assume that $pc_0$ large enough. If $u\in L^p_{\xi}(\mathbb{M}^N_{\varphi})$, then $u\equiv 0$ in $\mathbb{M}^N_{\varphi}$.
\item[(ii)] Assume that $b$ satisfies assumption \eqref{h3}. Then equation \eqref{problema} admits infinitely many classical solutions.
\end{itemize}

\end{corollary}

\subsection{Notations from Riemannian Geometry}

In this subsection we recall and fix some basic notations and useful definitions from Riemannian geometry that we will use in the rest of the paper (we refer e.g. to \cite{AMR,GHL,Grig}).

Throughout the paper, $M$ will denote an $N-$dimensional, complete, connected, noncompact, smooth Riemannian manifold without boundary and of infinite volume, endowed with a smooth Riemannian metric $g=\lbrace g_{ij}\rbrace$.  As usual, we consider the volume form, denoted by $d\mu$, induced by $g$ and we denote by $\mathrm{div}(X)$ the divergence of a smooth vector field $X$ on $M$, that is,  in local coordinates
\begin{equation*}
\mathrm{div} (X)=\dfrac{1}{\sqrt{\vert g\vert}}\sum_{i=1}^{N} \partial_i \left( \sqrt{\vert g\vert} X^i\right)\, ,
\end{equation*}
where $\vert g\vert=\mathrm{det}(g_{ij})(\geq 0)$. 

We also denote by $\nabla$ and $\Delta$ the Riemannian gradient and the Laplace-Beltrami operator on $M$ with respect to $g$, that is, in local coordinates,
$$
(\nabla u)^i=\sum_{j=1}^{N}g^{ij}\partial_j u \, ,
$$
and
$$
\Delta u=\mathrm{div}(\nabla u)=\dfrac{1}{\sqrt{\vert g\vert}}\sum_{i,j=1}^{N} \partial_i \left( \sqrt{\vert g\vert} g^{ij}\partial_j u\right)\, ,
$$
for any $C^2-$function $u:M\rightarrow\mathbb{R}$, where $\lbrace g^{ij}\rbrace$ denotes the inverse of the metric tensor $g$.

Due to the divergence structure of the Laplace-Beltrami operator we have the following integration by parts formula
\begin{equation}\label{parts}
\int_M v\Delta u\, d\mu=-\int_M \langle\nabla u,\nabla v\rangle\, d\mu\, ,
\end{equation}
for any $C^2-$functions $v,u:M\rightarrow\mathbb{R}$, with either $v$ or $u$ compactly supported and where $\langle\cdot,\cdot\rangle$ is the scalar product induced by $g$. 

Moreover, given a vector field $X$, we denote 
$$
\vert X\vert=\sqrt{\langle X,X\rangle}\,. 
$$
Observe that, for any $C^2-$functions $u:M\rightarrow\mathbb{R}$ and $\psi:\mathbb{R}\rightarrow\mathbb{R}$ we have
\begin{equation}\label{Laplacian_comp}
\Delta[\psi(u)]=\psi'(u)\Delta u + \psi''(u)\vert\nabla u\vert^{2}\,.
\end{equation}
Hence if $u$ is a classical solution of \eqref{problema} then 
\begin{equation}\label{soluzione_composta}
\Delta[\psi(u)]+ g(b,\nabla [\psi(u)]) - c\psi(u)=\psi''(u)\vert\nabla u\vert^2- c[\psi(u)-\psi'(u)u]\, . 
\end{equation}
Lastly, we recall the definition of model manifold (see e.g. \cite{GW,Petersen,PigS})

\begin{definition}\label{def_modello}
An $N-$dimensional model manifold with warping function $\varphi$ is a Riemannian manifold $\mathbb{M}^N_{\varphi}$ diffeomorphic to $\mathbb{R}^N$ and endowed with a rotationally symmetric Riemannian metric. The model $\mathbb{M}^N_{\varphi}$ is realized as the quotient space $[0,+\infty)\times\mathbb{S}^{N-1}/\backsim$, where $\backsim$ identifies $\lbrace 0\rbrace\times\mathbb{S}^{N-1}$ with the pole $o$ of the space, and the Riemannian metric has the expression
$$
g=dr^2+\varphi^2(r)d\theta^2\, 
$$
where $r(x)=\mathrm{dist}(x,o)$,  $d\theta^2$ denotes the standard metric of $\mathbb{S}^{N-1}$ and $\varphi:[0,+\infty)\rightarrow[0,+\infty)$ is a smooth function satisfying the following conditions: 
\begin{itemize}
\item $\varphi(r)>0$ for all $r>0$,
\item $\varphi'(0)=1$, 
\item $\varphi^{(2k)}(0)=0$ for every $k\geq 0$.
\end{itemize}
\end{definition}
Observe that the Euclidean space $\mathbb{R}^N$ and the (generalized) hyperbolic space $\mathbb{H}^N$ of constant curvature $-c<0$ can be seen as model manifolds (take $\varphi(r)=r$ and $\varphi(r)=\frac{\sinh(\sqrt{c}r)}{\sqrt{c}}$, respectively). Moreover, due to their structure the Laplace-Beltrami operator on $\mathbb{M}^N_{\varphi}$ can be written in the following way
\begin{equation}\label{Laplaciano_modelli}
\Delta = \dfrac{\partial^2}{\partial^2 r} + (N-1)\dfrac{\varphi'}{\varphi}\dfrac{\partial}{\partial r} + \dfrac{1}{\varphi^2}\Delta_\theta\, , 
\end{equation}
where $\Delta_\theta$ denotes the Laplace-Beltrami operator on $\mathbb{S}^{N-1}$. In particular, 
\begin{equation}\label{Laplaciano_dist}
\Delta r= (N-1)\dfrac{\varphi'}{\varphi}\, .
\end{equation}
Moreover, being the volume form on $\mathbb{M}^N_{\varphi}$ given by
$$
d\mu=\varphi^{N-1}(r)\, dr\, d\theta\, ,
$$
one has the following: 
\begin{equation}\label{volume_balle}
V(o,r)=c_N\int_{0}^{r}\varphi^{N-1}(t)\, dt \, ,
\end{equation}
where, $c_N$ is the area of the $(N-1)-$dimensional unit sphere and $o$ is the pole of $\mathbb{M}^N_{\varphi}$.

\

\section{Proof of Theorem \ref{teo1}}\label{p1}\setcounter{equation}{0}

\

Let $M$ be a complete, noncompact, Riemannian manifold such that \eqref{vol1} is satisfied. Firstly, we state a general criterion for uniqueness of solutions to equation \eqref{problema} in $L^p_{\psi}(M)$ for $\psi$ as in \eqref{eq12}. We suppose that there exists a positive function $\zeta\in C^1(M)$, which solves the following
\begin{equation}\label{eq40}
\delta |\nabla \zeta|^2-\langle b\,,\,\nabla \zeta\rangle-  \di b- c\,p < 0 \quad \textrm{in}\;\; M\,;
\end{equation}
for some $\delta>0$ sufficiently large and for some $p>1$. We mention that, $\delta=\delta(\varepsilon)$ where $\varepsilon>0$ will be specified in the proof of Theorem \ref{teo1}. Such inequality is meant in the sense of Definition \ref{defsole}-$(ii)$, where instead of $``="$ we have $``<"$.

\begin{proposition} \label{prop3}
Let $M$ be a complete, noncompact, Riemannian manifold such that \eqref{vol1} is satisfied. Moreover, assume that \eqref{h0} and \eqref{h2} hold. Let $u$ be a solution to equation \eqref{problema}. Assume that there exists a positive function $\zeta\in C^1(M)$ such that $\phi=e^{\zeta}$ satisfies 
\begin{equation}\label{eq42c}
\phi(x)\leq C \psi(x)\quad \textrm{for all}\;\; x\in M;
\end{equation}
for some constant $C>0$ and $\psi$ as in \eqref{eq12}.  Moreover, assume that $\zeta$ solves \eqref{eq40}. If $u\in L^p_\psi(M)$, then
$$
u\equiv 0\quad \textrm{in}\;\; M\,.
$$
\end{proposition}

\subsection{Proof of Proposition \ref{prop3}}

In order to show Proposition \ref{prop3} we need the following 

\begin{lemma}\label{lemma4}
Let $M$ be a complete, noncompact, Riemannian manifold. Assume that \eqref{h0} and \eqref{h2} hold. Let $\phi\in C^1(M)$, $\phi>0$ be such that \eqref{eq42c} holds.
Let $v\in L^1_{\psi}(M)$. Then
\begin{equation}\label{eq41bis}
\lim_{R\to\infty}\int_{M}|v(x)| \phi|\nabla\gamma_R(x)|^2\,d\mu = 0, 
\end{equation}
and
\begin{equation}\label{eq41tris}
\lim_{R\to\infty}\left[\int_{D_+}|v(x)| \phi\,\gamma_R(x)\left\langle b(x)\,,\,\nabla\gamma_R(x)\right\rangle\,d\mu \right]= 0;
\end{equation}
 for every cut-off function $\gamma_R\in C_c^{\infty}(B_R)$ and where $D_+$ defined as in \eqref{d}.
\end{lemma}

\begin{remark}\label{rem31}
Observe that the quantity into the brackets of formula \eqref{eq41tris} is negative due to \eqref{d}.
\end{remark}

\begin{proof}[Proof of Lemma \ref{lemma4}]
Let us consider a cut-off function $\gamma\in C^{\infty}([0,+\infty))$, $0\le \gamma\le1$ such that
$$ \gamma(r):=\begin{cases} 1\quad&0\le r\le \frac 12\\ 0\quad &r\ge 1\,,\end{cases}$$
where $r$ is the geodesic distance given by \eqref{r}. Then, let us define
\begin{equation}\label{gam1}
\gamma_R(x):=\gamma\left(\frac rR\right)\quad \text{for all}\,\,x\in M.
\end{equation}
To show \eqref{eq41bis}, let us observe that, for any $x\in B_R\setminus B_{R/2}$, for some $\tilde C>0$
\begin{equation}\label{eq33a}
|\nabla\gamma_R|^2=\left[\gamma'_R\left(\frac rR\right)\right]^2\frac 1R|\nabla r|^2\le\,\frac {\tilde C}{R^2}\le \frac{\tilde C}{r^2}.
\end{equation}
Therefore, by combining \eqref{eq42c} with \eqref{eq33a}, we obtain
$$
\begin{aligned}
\int_M|v|\phi\,|\nabla\gamma_R|^2\,d\mu\,&\le\, \int_{B_R\setminus B_{R/2}}|v|\phi\,\frac{\tilde C}{r^2}\,d\mu\le\,\tilde C C\frac{4}{R^2}\int_{B_R\setminus B_{R/2}}|v|\psi\,d\mu\,.
\end{aligned}
$$
Then, since $v\in L^1_{\psi}(M)$, from the latter we obtain \eqref{eq41bis} by letting $R\to\infty$. On the other hand, to show \eqref{eq41tris}, let us observe that, due to \eqref{h0}, for any $x\in D_+\cap (B_R\setminus B_{R/2})$, for some $\bar C>0$
\begin{equation}\label{eq33b}
-\left\langle b(x),\nabla \gamma_R(x)\right\rangle= -\gamma'\left(\frac{r}{R}\right)\frac{1}{R}\,\left\langle b(x),\nabla r\right\rangle\,\le \bar C K,
\end{equation}
where, $D_+$ has been defined in \eqref{d}. Then, due to \eqref{eq42c} together with \eqref{eq33b}, we get
\begin{equation}\label{eq33}
\begin{aligned}
-\int_{D_+} &|v(x)|\phi(x)\,\gamma_R\langle b(x),\nabla \gamma_R\rangle\,d\mu  \\
&\le \int_{D_+\cap (B_R\setminus B_{R/2})} |v(x)|\,\phi(x)\,\bar CK\,d\mu \le \bar C CK \int_{D_+\cap (B_R\setminus B_{R/2})} |v(x)| \psi(x)\,d\mu.
\end{aligned}
\end{equation}
Then, since $v\in L^1_{\psi}(M)$, we obtain from \eqref{eq33} that
\begin{equation*}\label{eq34}
\lim_{R\to\infty}\left[- \int_{D_+} |v(x)|\,\phi(x)\,\gamma_R\,\left\langle b(x),\nabla \gamma_R(x)\right\rangle\,d\mu\right] =0.
\end{equation*}
This completes the proof of the Lemma.
\end{proof}

\begin{proof}[Proof of Proposition \ref{prop3}]
For any $\alpha>0$ small enough, we define
\begin{equation}\label{G}
G_{\alpha}(u):=(u^2+\alpha)^{\frac p4}\,.
\end{equation}
Then we take $v\in C^1(M)$ such that
$$
v:=G_{\alpha}(u)\gamma_R^2\,e^{\zeta},
$$
where $\gamma_R$ has been defined in \eqref{gam1}. Now, we consider left-hand side of equation \eqref{problema} with $u=G_{\alpha}(u)$, i.e.
\begin{equation}\label{eq42a}
\Delta[G_{\alpha}(u)] + \left\langle b\,,\,\nabla [G_{\alpha}(u)] \right\rangle - c \,G_{\alpha}(u)\,.
\end{equation}
Then, we multiply it by $v$, we integrate over $M$ and we apply the integration by parts formula \eqref{parts}. We obtain
\begin{equation}\label{eq42}
\begin{aligned}
\int_M &\left\{\Delta[G_{\alpha}(u)] + \left\langle b\,,\,\nabla [G_{\alpha}(u)] \right\rangle - c \,G_{\alpha}(u)\right\} \,v\,d\mu \\
& =-\int_M |\nabla G_{\alpha}(u)|^2\gamma_R^2e^{\zeta}\,d\mu -\int_M 2\left\langle \nabla G_{\alpha}(u),\nabla \gamma_R\right\rangle G_{\alpha}(u)\gamma_R e^{\zeta}\,d\mu\\
& \quad-\int_M\left\langle \nabla G_{\alpha}(u),\nabla \zeta \right\rangle G_{\alpha}(u)\gamma_R^2e^{\zeta}\,d\mu \\
&\quad+\frac 12 \int_M \left\langle b\,,\,\nabla [G_{\alpha}(u)] \right\rangle G_{\alpha}(u)\gamma_R^2e^{\zeta}\,d\mu - \frac 12 \int_M \left\langle b\,,\,\nabla [G_{\alpha}(u)] \right\rangle G_{\alpha}(u)\gamma_R^2e^{\zeta}\, d\mu\\
&\quad-\frac 12 \int_M 2\left\langle b\,,\,\nabla \gamma_R \right\rangle G^2_{\alpha}(u)\gamma_Re^{\zeta}\,d\mu - \frac 12 \int_M \left\langle b\,,\,\nabla \zeta\right\rangle G^2_{\alpha}(u)\gamma_R^2e^{\zeta}\, d\mu\\
&\quad-\frac 12 \int_M\operatorname{div} b \,G_{\alpha}^2(u)\gamma_R^2e^{\zeta}\,d\mu - \int_M c\,G_{\alpha}^2(u)\gamma_R^2e^{\zeta}\,d\mu.
\end{aligned}
\end{equation}
By using Young's inequality, for every $\varepsilon>0$, \eqref{eq42} reduces to
\begin{equation}\label{eq42f}
\begin{aligned}
\int_M &\left\{\Delta[G_{\alpha}(u)] + \left\langle b\,,\,\nabla [G_{\alpha}(u)] \right\rangle - c \,G_{\alpha}(u)\right\} \,v\,d\mu \\
&\le-\varepsilon\int_M |\nabla G_{\alpha}(u)|^2\gamma_R^2e^{\zeta}\,d\mu -(1-\varepsilon)\int_M |\nabla G_{\alpha}(u)|^2\gamma_R^2e^{\zeta}\,d\mu \\
&\quad + \frac{\varepsilon}{2}\int_M| \nabla G_{\alpha}(u)|^2\gamma_R^2 e^{\zeta}\,d\mu+\frac{2}{\varepsilon}\int_M|\nabla \gamma_R|^2 G_{\alpha}(u)^2 e^{\zeta}\,d\mu\\
& \quad+\frac{\varepsilon}{2}\int_M| \nabla G_{\alpha}(u)|^2 \gamma_R^2 e^{\zeta}\,d\mu+\frac{1}{2\varepsilon}\int_M |\nabla\zeta|^2 G_{\alpha}(u)^2 \gamma_R^2 e^{\zeta}\,d\mu \\
&\quad+\frac 12 \int_M \left\langle b\,,\,\nabla [G_{\alpha}(u)] \right\rangle G_{\alpha}(u)\gamma_R^2e^{\zeta}\,d\mu - \frac 12 \int_M \left\langle b\,,\,\nabla [G_{\alpha}(u)] \right\rangle G_{\alpha}(u)\gamma_R^2e^{\zeta}\, d\mu\\
&\quad-\frac 12 \int_M 2\left\langle b\,,\,\nabla \gamma_R \right\rangle G^2_{\alpha}(u)\gamma_Re^{\zeta}\,d\mu - \frac 12 \int_M \left\langle b\,,\,\nabla \zeta\right\rangle G^2_{\alpha}(u)\gamma_R^2e^{\zeta}\, d\mu\\
&\quad-\frac 12 \int_M\operatorname{div} b \,G_{\alpha}^2(u)\gamma_R^2e^{\zeta}\,d\mu - \int_M c\,G_{\alpha}^2(u)\gamma_R^2e^{\zeta}\,d\mu\\
&\le -(1-\varepsilon)\int_M|G'_{\alpha}(u)|^2 |\nabla u|^2\gamma_R^2e^{\zeta}\,d\mu+ \frac{2}{\varepsilon} \int_M  |\nabla \gamma_R|^2 G^2_{\alpha}(u) e^{\zeta}\,d\mu \\
&\quad+ \frac 1{2\varepsilon} \int_M |\nabla \zeta|^2 G^2_{\alpha}(u)\gamma_R^2e^{\zeta}\,d\mu - \int_M \left\langle b\,,\,\nabla \gamma_R \right\rangle G^2_{\alpha}(u)\gamma_Re^{\zeta}\,d\mu\\
&\quad  - \frac 12 \int_M \left\langle b\,,\,\nabla \zeta\right\rangle G^2_{\alpha}(u)\gamma_R^2e^{\zeta}\, d\mu-\frac 12 \int_M\operatorname{div} b \,G_{\alpha}^2(u)\gamma_R^2e^{\zeta}\,d\mu \\
&\quad - \int_M c\,G_{\alpha}^2(u)\gamma_R^2e^{\zeta}\,d\mu\,.\\
\end{aligned}
\end{equation}
On the other hand, formula \eqref{eq42a} can be reduced as follows
\begin{equation}\label{eq43b}
\begin{aligned}
\Delta&[G_{\alpha}(u)] +\langle b,\nabla [G_{\alpha}(u)]\rangle - c \,G_{\alpha}(u)\\
&=G_{\alpha}'(u)\Delta u+G_{\alpha}''(u)|\nabla u|^2+G_{\alpha}'(u)\langle b,\nabla u\rangle-c\,G_{\alpha}(u)+c\,G_{\alpha}'(u)u-c\,G_{\alpha}'(u)u\\
&=G_{\alpha}'(u)\left[\Delta u+\langle b,\nabla u\rangle-cu\right] +G_{\alpha}''(u)|\nabla u|^2-c\,G_{\alpha}(u)+c\,G_{\alpha}'(u)u\\
&= G_{\alpha}''(u)|\nabla u|^2+G_{\alpha}'(u)\,c\,u-c\,G_{\alpha}(u)\,,
\end{aligned}
\end{equation}
where we used \eqref{Laplacian_comp}, \eqref{soluzione_composta} and \eqref{problema}. From \eqref{eq42f} and \eqref{eq43b}, 
\begin{equation}\label{eq42g}
\begin{aligned}
\int_M &\left[G_{\alpha}''(u)|\nabla u|^2+G_{\alpha}'(u)\,c\,u-c\,G_{\alpha}(u)\right] \,G_{\alpha}(u)\gamma_R^2e^{\zeta}\,d\mu\\
&\le -(1-\varepsilon)\int_M|G'_{\alpha}(u)|^2 |\nabla u|^2\gamma_R^2e^{\zeta}\,d\mu+ \frac{2}{\varepsilon} \int_M  |\nabla \gamma_R|^2 G^2_{\alpha}(u) e^{\zeta}\,d\mu \\
&\quad+ \frac 1{2\varepsilon} \int_M |\nabla \zeta|^2 G^2_{\alpha}(u)\gamma_R^2e^{\zeta}\,d\mu - \int_M \left\langle b\,,\,\nabla \gamma_R \right\rangle G^2_{\alpha}(u)\gamma_Re^{\zeta}\,d\mu\\
&\quad  - \frac 12 \int_M \left\langle b\,,\,\nabla \zeta\right\rangle G^2_{\alpha}(u)\gamma_R^2e^{\zeta}\, d\mu-\frac 12 \int_M\operatorname{div} b \,G_{\alpha}^2(u)\gamma_R^2e^{\zeta}\,d\mu \\
&\quad - \int_M c\,G_{\alpha}^2(u)\gamma_R^2e^{\zeta}\,d\mu\,.\\
\end{aligned}
\end{equation}
By rearranging the terms in \eqref{eq42g}, we get
\begin{equation}\label{eq43}
\begin{aligned}
\int_M &\left[G_{\alpha}'(u)\,c\,u-c\,G_{\alpha}(u)\right] \,G_{\alpha}(u)\gamma_R^2e^{\zeta}\,d\mu\\
& \le -\int_M \left\{(1-\varepsilon)|G'_{\alpha}(u)|^2+G''_{\alpha}(u)G_{\alpha}(u)\right\}|\nabla u|^2\gamma_R^2e^{\zeta}\,d\mu \\
& \quad + \frac{2}{\varepsilon} \int_M  |\nabla \gamma_R|^2 G^2_{\alpha}(u) e^{\zeta}\,d\mu + \frac 1{2\varepsilon} \int_M |\nabla \zeta|^2 G^2_{\alpha}(u)\gamma_R^2e^{\zeta}\,d\mu\\
&\quad - \int_M \left\langle b\,,\,\nabla \gamma_R \right\rangle G^2_{\alpha}(u)\gamma_Re^{\zeta}\,d\mu - \frac 12 \int_M \left\langle b\,,\,\nabla \zeta\right\rangle G^2_{\alpha}(u)\gamma_R^2e^{\zeta}\, d\mu\\
&\quad-\frac 12 \int_M\operatorname{div} b \,G_{\alpha}^2(u)\gamma_R^2e^{\zeta}\,d\mu - \int_M c\,G_{\alpha}^2(u)\gamma_R^2e^{\zeta}\,d\mu\,.\\
\end{aligned}
\end{equation}
On the left-hand side of \eqref{eq43}, by substituting the definition of $G_{\alpha}$ given in \eqref{G}, we have
\begin{equation}\label{eq43f}
\begin{aligned}
\int_M &\left[G_{\alpha}'(u)\,c\,u-c\,G_{\alpha}(u)\right] \,G_{\alpha}(u)\gamma_R^2e^{\zeta}\,d\mu\\
&=\int_M (u^2+\alpha)^{\frac p4-1}\,c\left[\left(\frac p2-1\right)u^2-\alpha\right] (u^2+\alpha)^{\frac p4}\gamma_R^2\,e^{\zeta}\,d\mu\,.
\end{aligned}
\end{equation}
Similarly, on the right-hand side of \eqref{eq43f}, by exploiting the definition of $G_{\alpha}$ given in \eqref{G} and by choosing $\varepsilon>0$ sufficiently small,  for any $p>1$, we have that
\begin{equation}\label{eq43g}
-\int_M \left\{(1-\varepsilon)|G'_{\alpha}(u)|^2+G''_{\alpha}(u)G_{\alpha}(u)\right\}|\nabla u|^2\gamma_R^2e^{\zeta}\,d\mu\le 0,
\end{equation}
since, by choosing $0<\varepsilon\le 2-\frac 2p$,
$$
\begin{aligned}
(1-\varepsilon)|G'_{\alpha}(u)&|^2+G''_{\alpha}(u)G_{\alpha}(u)\\
&=(1-\varepsilon)\frac{p^2}{4}(u^2+\alpha)^{\frac p2-2}u^2+\frac p2 (u^2+\alpha)^{\frac p2-1}+ p\left(\frac p4 -1\right)(u^2+\alpha)^{\frac p2-2} u^2\\
&=(2-\varepsilon)\frac{p^2}{4}(u^2+\alpha)^{\frac p2-2}u^2+\frac p2 (u^2+\alpha)^{\frac p2-1}-p(u^2+\alpha)^{\frac p2-2} u^2\\
&=(u^2+\alpha)^{\frac p2-2}\left[u^2\left((2-\varepsilon)\frac{p^2}{4}-\frac p2\right)+\frac p2\alpha\right]\\
&\ge 0\,.
\end{aligned}
$$
Hence, \eqref{eq43}, combined with \eqref{eq43f} and \eqref{eq43g}, reduces to
\begin{equation}\label{eq43c}
\begin{aligned}
\int_M (u^2+\alpha)^{\frac p4-1}\,&c\left[\left(\frac p2-1\right)u^2-\alpha\right] (u^2+\alpha)^{\frac p4}\gamma_R^2\,e^{\zeta}\,d\mu\\
& \le+ \frac{2}{\varepsilon} \int_M  |\nabla \gamma_R|^2 G^2_{\alpha}(u) e^{\zeta}\,d\mu + \frac 1{2\varepsilon} \int_M |\nabla \zeta|^2 G^2_{\alpha}(u)\gamma_R^2e^{\zeta}\,d\mu\\
&\quad - \int_M \left\langle b\,,\,\nabla \gamma_R \right\rangle G^2_{\alpha}(u)\gamma_Re^{\zeta}\,d\mu - \frac 12 \int_M \left\langle b\,,\,\nabla \zeta\right\rangle G^2_{\alpha}(u)\gamma_R^2e^{\zeta}\, d\mu\\
&\quad-\frac 12 \int_M\operatorname{div} b \,G_{\alpha}^2(u)\gamma_R^2e^{\zeta}\,d\mu - \int_M c\,G_{\alpha}^2(u)\gamma_R^2e^{\zeta}\,d\mu\\
\end{aligned}
\end{equation}
Letting $\a\to 0^+$ 
$$
(u^2+\alpha)^{\frac p4-1}\,c\left[\left(\frac p2-1\right)u^2-\alpha\right] (u^2+\alpha)^{\frac p4}\longrightarrow\left(\frac p2 -1\right)|u|^p\,,
$$
and
$$
G_{\alpha}(u)=(u^2+\alpha)^{\frac p4}\longrightarrow |u|^{\frac p2}\,.
$$
Hence, by the dominated convergence theorem, \eqref{eq43c} becomes
\begin{equation*}
\begin{aligned}
\left(\frac p2 -1\right)\int_M c\,|u|^p\,\gamma_R^2\,e^{\zeta}\,d\mu&\le+ \frac{2}{\varepsilon} \int_M  |\nabla \gamma_R|^2 |u|^p e^{\zeta}\,d\mu + \frac 1{2\varepsilon} \int_M |\nabla \zeta|^2 |u|^p\gamma_R^2e^{\zeta}\,d\mu\\
&\quad - \int_M \left\langle b\,,\,\nabla \gamma_R \right\rangle |u|^p\gamma_Re^{\zeta}\,d\mu - \frac 12 \int_M \left\langle b\,,\,\nabla \zeta\right\rangle |u|^p\gamma_R^2e^{\zeta}\, d\mu\\
&\quad-\frac 12 \int_M\operatorname{div} b \,|u|^p\gamma_R^2e^{\zeta}\,d\mu - \int_M c\,|u|^p\gamma_R^2e^{\zeta}\,d\mu,\\
\end{aligned}
\end{equation*}
which, defining $\delta:=\frac 1\varepsilon$, is equivalent to
\begin{equation}\label{eq44}
\begin{aligned}
\frac 12\int_M  |u|^p\gamma_R^2e^{\zeta} &\left[-\delta |\nabla \zeta|^2+\left\langle b\,,\,\nabla \zeta\right\rangle+  \di b+ c\, p\right]\, d\mu\\
&\le \int_M |u|^pe^{\zeta}\left[2\delta |\nabla \gamma_R|^2- \left\langle b\,,\,\nabla \gamma_R \right\rangle\gamma_R\right]\,d\mu.
\end{aligned}
\end{equation}
Due to \eqref{eq44}, and the definitions of $D_+$ and $D_-$ in \eqref{d}, we obtain
\begin{equation}\label{eq46}
\begin{aligned}
& \frac 12\int_M |u|^p\gamma_R^2\,e^{\zeta} \left[-\delta |\nabla \zeta|^2+\left\langle b\,,\,\nabla \zeta\right\rangle+  \di b+ c\,p\right]\, d\mu\\
&\quad \le \int_{M} 2\delta\,|u|^p e^{\zeta} |\nabla \gamma_R|^2\,d\mu-\int_{D_+}|u|^p e^{\zeta} \left\langle b,\nabla \gamma_R \right\rangle\gamma_R\,d\mu -\int_{D_-}|u|^p e^{\zeta} \left\langle b,\nabla \gamma_R \right\rangle\gamma_R\,d\mu.
\end{aligned}
\end{equation}
Observing that, for all $x\in D_-$
$$
- \langle b(x)\,,\,\nabla \gamma_R(x)\rangle \,\le 0,
$$
we can rewrite \eqref{eq46} as follows
\begin{equation}\label{eq46a}
\begin{aligned}
\frac 12\int_M  |u|^p\gamma_R^2\,e^{\zeta} &\left[-\delta |\nabla \zeta|^2+\left\langle b\,,\,\nabla \zeta\right\rangle+  \di b+ c\,p\right]\, d\mu\\
&\quad \le \int_{M} 2\delta\,|u|^p e^{\zeta} |\nabla \gamma_R|^2\,d\mu-\int_{D_+}|u|^p e^{\zeta} \left\langle b,\nabla \gamma_R \right\rangle\gamma_R\,d\mu.
\end{aligned}
\end{equation}
We now set $\phi(x):=e^{\zeta(x)}$ for any $x\in M$. Observe that, $\phi$ satisfies assumption \eqref{eq42c}. Hence, we can apply Lemma \ref{lemma4} with this choice of $\phi$. From Lemma \ref{lemma4} and the monotone convergence theorem, sending $R\to \infty$ in \eqref{eq46a} we get
\begin{equation}\label{eq47}
 \frac 12 \int_M |u|^p\gamma_R^2\,e^{\zeta} \left[-\delta |\nabla \zeta|^2+\left\langle b\,,\,\nabla \zeta\right\rangle+  \di b+ c\,p\right]\, d\mu\leq 0\,.
\end{equation}
From \eqref{eq47} and \eqref{eq40}, since $|u|^p\geq 0$, we can infer that $u\equiv 0$ in $M$. This completes the proof.
\end{proof}

We are now ready to prove Theorem \ref{teo1}, due to Proposition \ref{prop3} it suffices to construct a subsolution to \eqref{eq40}.

\begin{proof}[Proof of Theorem \ref{teo1}]
Let $\phi=\phi(r):=\psi(r)$ for any $r>0$, where $\psi$ has been defined in \eqref{eq12}, with $\beta>\alpha$. Moreover, we set $\zeta(r):=-\beta r$. Hence, $\phi(r)\equiv e^{\zeta(r)}$.
At first observe that $\phi$ satisfies \eqref{eq42c}. Moreover, $\zeta$ solves \eqref{eq40}, for properly chosen $\beta>\alpha$. Indeed, due to \eqref{h0} and \eqref{h2}, we have
\begin{equation}\label{eq45}
\begin{aligned}
\delta& |\nabla \zeta|^2-\left\langle b\,,\,\nabla \zeta\right\rangle-  \di b- c\,p\\
&=\delta \zeta'(r)^2|\nabla r|^2- \zeta'(r) \left\langle b\,,\,\nabla r\right\rangle -  \di b- c\,p\\
&\le -\left( \b^2\delta-\b K_1 -K +  p\,c_0\right).
\end{aligned}
\end{equation}
Then, by \eqref{eq45}, one has
$$
\delta |\nabla \zeta|^2-\left\langle b\,,\,\nabla \zeta\right\rangle-  \di b- c\,p<0,
$$ 
provided that
$$
p\,c_0>\b^2\delta+\b K.
$$
Thus, by Proposition \ref{prop3} the conclusion follows. 
\end{proof}

\

\section{Proof of Theorem \ref{teo3}}\label{p4}\setcounter{equation}{0}

\

Similarly to Theorem \ref{teo1}, we want to state a general criteria for uniqueness by means of assumption \eqref{h4}. Therefore, we state the next

\begin{proposition} \label{prop8}
Let $M$ be a complete, noncompact, Riemannian manifold such that \eqref{vol3} is satisfied. Moreover, assume that \eqref{h4} and \eqref{h2} hold. Let $u$ be a solution to equation \eqref{problema} and $\eta$ be as in \eqref{eq2222}. Assume that there exists a positive function $\zeta\in C^1(M)$ such that $\phi=e^{\zeta}$ satisfies 
\begin{equation}\label{eq82d}
\frac{\phi^2(x)}{1+r^2}\leq C \eta(x)\quad \textrm{for all}\;\; x\in M;
\end{equation}
and
\begin{equation}\label{eq82c}
\frac{\left\langle b(x),\nabla r\right\rangle}{(1+r)^{\theta}}\phi^2(x)\leq C \eta(x)\quad \textrm{for all}\;\; x\in D_+.
\end{equation}
Moreover, assume that $\zeta$ solves \eqref{eq40}. If $u\in L^p_\eta(M)$, then
$$
u\equiv 0\quad \textrm{in}\;\; M\,.
$$
\end{proposition}

Analogously to Lemma \ref{lemma4} the next lemma can be shown. The proof is similar and for this reason we only sketch it.

\begin{lemma}\label{lemma8}
Let $M$ be a complete, noncompact, Riemannian manifold such that \eqref{vol3} holds. Assume \eqref{h4} and \eqref{h2}. Let $\phi\in C^1(M)$, $\phi>0$ be such that \eqref{eq82d} and \eqref{eq82c} hold.
Let $v\in L^1_{\eta}(M)$. Then
\begin{equation}\label{eq81bis}
\lim_{R\to\infty}\int_{M}|v(x)| \phi^2|\nabla\gamma_R(x)|^2\,d\mu = 0, 
\end{equation}
and
\begin{equation}\label{eq81tris}
\lim_{R\to\infty}\left[\int_{D_+}|v(x)| \phi^2\,\gamma_R(x)\left\langle b(x)\,,\,\nabla\gamma_R(x)\right\rangle\,d\mu \right]= 0;
\end{equation}
 for every cut-off function $\gamma_R\in C_c^{\infty}(B_R)$ and where $D_+$ defined as in \eqref{d}.
\end{lemma}

\begin{proof}
Equality \eqref{eq81bis} can be obtained as \eqref{eq41bis} by using assumption \eqref{eq82d}. To show \eqref{eq81tris} it is sufficient to argue as in the proof of Lemma \ref{lemma4} by replacing assumptions \eqref{h0} and \eqref{eq42c} with \eqref{h4} and \eqref{eq82c} respectively. 
\end{proof}

\begin{proof}[Proof of Proposition \ref{prop8}]
To prove Proposition \ref{prop8} it is sufficient to argue as in the proof of Proposition \ref{prop3} by using Lemma \ref{lemma8} instead of Lemma \ref{lemma4}.

\end{proof}

\begin{proof}[Proof of Theorem \ref{teo3}]
Let $\phi=\phi(r):=\eta(r)$ for any $r>0$, where $\eta$ has been defined in \eqref{eq2222}. Moreover, we set $\zeta(r):=-\beta r^{\theta}$, $0<\theta<1$ as in \eqref{h4}. Hence, $\phi(r)\equiv e^{\zeta(r)}$.
At first observe that $\phi$ satisfies \eqref{eq82d} and \eqref{eq82c}. Moreover, $\zeta$ solves \eqref{eq40}, for properly chosen $\beta>\alpha$. Indeed, due to \eqref{h4} and \eqref{h2}, we have
\begin{equation}\label{eq45bis}
\begin{aligned}
\delta& |\nabla \zeta|^2-\left\langle b\,,\,\nabla \zeta\right\rangle-  \di b- c\,p\\
&=\delta \zeta'(r)^2|\nabla r|^2- \zeta'(r) \left\langle b\,,\,\nabla r\right\rangle -  \di b- c\,p\\
&\le \beta^2\theta^2\delta\,r^{2\theta-2}+\beta\theta\,K_1\,r^{\theta-1-\sigma}+K\,r^{\sigma-1}-c_0\,p\\
&\le -\left(-\beta^2\theta^2\delta-\b K_1\theta -K +  p\,c_0\right).
\end{aligned}
\end{equation}
Then, by \eqref{eq45bis}, one has
$$
\delta |\nabla \zeta|^2-\left\langle b\,,\,\nabla \zeta\right\rangle-  \di b- c\,p<0,
$$ 
provided that
$$
p\,c_0>\b^2\theta^2\delta+\b K_1\theta+\,K.
$$
Thus, by Proposition \ref{prop8} the conclusion follows. 
\end{proof}

\

\section{Proof of Theorem \ref{teo2}}\label{p2}\setcounter{equation}{0}

\

Let $M$ be a complete, noncompact, Riemannian manifold such that \eqref{vol2} is satisfied. Firstly, we state a general criterion for uniqueness of solutions to equation \eqref{problema} in $L^p_{\xi}(M)$ for $\xi$ as in \eqref{eq13}. We suppose that there exists a positive function $\phi\in C^1(M)$, which solves, the following 
\begin{equation}\label{eq50}
\hat\delta|\nabla\phi|^2 - \left \langle b,\,\nabla \phi \right \rangle \phi - \frac 12(\di b+ p\,c )\,\phi^2 < 0 \quad \textrm{in}\;\; M\,;
\end{equation}
for some $\hat\delta>0$ sufficiently large and for some $p>1$. We mention that, $\hat\delta=\hat\delta(\varepsilon)$ wher $\varepsilon>0$ will be specified in the proof of Theorem \ref{teo2}. Such inequality is meant in the sense of Definition \ref{defsole}-$(ii)$, where instead of $``="$ we have $``<"$.

\begin{proposition} \label{prop5}
Let $M$ be a complete, noncompact, Riemannian manifold such that \eqref{vol2} is satisfied. Moreover, assume that \eqref{h1} and \eqref{h2} hold. Let $u$ be a solution to equation \eqref{problema} and $\xi$ be as in \eqref{eq13}. Assume that there exists a positive function $\phi\in C^1(M)$, which solves \eqref{eq50} and such that
\begin{equation}\label{eq52d}
\frac{\phi^2(x)}{1+r^2}\leq C \xi(x)\quad \textrm{for all}\;\; x\in M;
\end{equation}
and
\begin{equation}\label{eq52c}
\frac{\left\langle b(x),\nabla r\right\rangle}{1+r}\phi^2(x)\leq C \xi(x)\quad \textrm{for all}\;\; x\in D_+
\end{equation}
If $u\in L^p_\xi(M)$, then
$$
u\equiv 0\quad \textrm{in}\;\; M\,.
$$
\end{proposition}

\subsection{Proof of Proposition \ref{prop5}}
Analogously to Lemma \ref{lemma4} and Lemma \ref{lemma8} the next lemma can be shown.
\begin{lemma}\label{lemma5}
Let $M$ be a complete, noncompact, Riemannian manifold such that \eqref{vol2} holds. Assume \eqref{h1} and \eqref{h2}. Let $\phi\in C^1(M)$, $\phi>0$ be such that \eqref{eq52d} and \eqref{eq52c} hold.
Let $v\in L^1_{\xi}(M)$. Then
\begin{equation}\label{eq51bis}
\lim_{R\to\infty}\int_{M}|v(x)| \phi^2|\nabla\gamma_R(x)|^2\,d\mu = 0, 
\end{equation}
and
\begin{equation}\label{eq51tris}
\lim_{R\to\infty}\left[\int_{D_+}|v(x)| \phi^2\,\gamma_R(x)\left\langle b(x)\,,\,\nabla\gamma_R(x)\right\rangle\,d\mu \right]= 0;
\end{equation}
 for every cut-off function $\gamma_R\in C_c^{\infty}(B_R)$ and where $D_+$ defined as in \eqref{d}.
\end{lemma}

\begin{proof}
Equality \eqref{eq51bis} can be obtained as \eqref{eq41bis} by using assumption \eqref{eq51bis}. To show \eqref{eq51tris} it is sufficient to argue as in the proof of Lemma \ref{lemma4} by replacing assumptions \eqref{h0} and \eqref{eq42c} with \eqref{h1} and \eqref{eq52c} respectively. 
\end{proof}

\begin{proof}[Proof of Proposition \ref{prop5}]
For any $\alpha>0$ small enough, Let $G_{\alpha}$ be defined as in \eqref{G}. Then we take $v\in C^1(M)$ such that
$$
v:=G_{\alpha}(u)\gamma_R^2\,\phi^2, \quad \text{for all}\,\,\,x\in M
$$ 
where $\gamma_R$ has been defined in \eqref{gam1}. Now, we consider the left-hand side of equation \eqref{problema} with $u=G_{\alpha}(u)$, i.e.
\begin{equation}\label{eq52a}
\Delta[G_{\alpha}(u)] + \left\langle b\,,\,\nabla G_{\alpha}(u) \right\rangle - c \,G_{\alpha}(u)\,.
\end{equation}
Then, we multiply it by $v$, we integrate over $M$ and we use the integration by parts formula \eqref{parts}
\begin{equation*}
\begin{aligned}
\int_M &\left\{\Delta[G_{\alpha}(u)] + \left\langle b\,,\,\nabla G_{\alpha}(u) \right\rangle - c \,G_{\alpha}(u)\right\} \,v\,d\mu \\
&=-\int_M \left\langle \nabla G_{\alpha},\,\nabla v\right\rangle\,d\mu+\frac 12 \int_M \left\langle b\,,\,\nabla G_{\alpha}(u) \right\rangle\, v\,d\mu  -\frac 12 \int_M \left\langle b\,,\,\nabla v \right\rangle G_{\alpha}(u) \\& \quad - \int_M \left(\frac 12 \operatorname{div} b+c\right )\, v\,d\mu\\
& =-\varepsilon\int_M |\nabla G_{\alpha}(u)|^2\gamma_R^2\phi^2\,d\mu -(1-\varepsilon)\int_M |\nabla G_{\alpha}(u)|^2\gamma_R^2\phi^2\,d\mu \\
&\quad - 2\int_M\langle\nabla G_{\alpha}(u),\nabla \gamma_R\rangle G_{\alpha}(u)\gamma_R \phi^2\,d\mu-2\int_M\langle \nabla G_{\alpha}(u),\nabla\phi\rangle G_{\alpha}(u) \gamma_R^2  \phi\,d\mu\\
&\quad+\frac 12 \int_M \left\langle b\,,\,\nabla [G_{\alpha}(u)] \right\rangle G_{\alpha}(u)\gamma_R^2\phi^2\,d\mu - \frac 12 \int_M \left\langle b\,,\,\nabla [G_{\alpha}(u)] \right\rangle G_{\alpha}(u)\gamma_R^2\phi^2\, d\mu\\
&\quad-\frac 12 \int_M 2\left\langle b\,,\,\nabla \gamma_R \right\rangle G^2_{\alpha}(u)\gamma_R\phi^2\,d\mu - \frac 12 \int_M 2 \left\langle b\,,\,\nabla \phi\right\rangle G^2_{\alpha}(u)\gamma_R^2\phi^2\, d\mu\\
&\quad-\frac 12 \int_M\operatorname{div} b \,G_{\alpha}^2(u)\gamma_R^2\phi^2\,d\mu - \int_M c\,G_{\alpha}^2(u)\gamma_R^2\phi^2\,d\mu \, . 
\end{aligned}
\end{equation*}
By Young's inequality we obtain
\begin{equation}\label{eq52}
\begin{aligned}
\int_M &\left\{\Delta[G_{\alpha}(u)] + \left\langle b\,,\,\nabla G_{\alpha}(u) \right\rangle - c \,G_{\alpha}(u)\right\} \,v\,d\mu \\
&\le \left(-\varepsilon+\frac {\varepsilon}2 +\frac{\varepsilon}2 \right) \int_M |\nabla G_{\alpha}(u)|^2\gamma_R^2\,\phi^2\,d\mu-\left(1-\varepsilon \right) \int_M |\nabla G_{\alpha}(u)|^2\gamma_R^2\,\phi^2\,d\mu\\
&\quad +\frac2{\varepsilon}  \int_M |\nabla \gamma_R|^2 G^2_{\alpha}(u)\phi^2\,d\mu+ \frac2{\varepsilon}   \int_M |\nabla \phi|^2 G^2_{\alpha}(u)\gamma_R^2\,d\mu\\
&\quad - \int_M \left\langle b\,,\,\nabla \gamma_R \right\rangle G^2_{\alpha}(u)\gamma_R\,\phi^2\,d\mu -  \int_M \left\langle b\,,\,\nabla \phi\right\rangle G^2_{\alpha}(u)\gamma_R^2\,\phi\, d\mu\\
&\quad-\frac 12 \int_M\operatorname{div} b \,G_{\alpha}^2(u)\gamma_R^2\phi^2\,d\mu - \int_M c\,G_{\alpha}^2(u)\gamma_R^2\phi^2\,d\mu\, .
\end{aligned}
\end{equation}
On the other hand, arguing as in \eqref{eq43b}, \eqref{eq52a} can be reduced as follows
\begin{equation}\label{eq53b}
\begin{aligned}
\Delta[G_{\alpha}(u)] &+\langle b,\nabla [G_{\alpha}(u)]\rangle - c \,G_{\alpha}(u)\\
&=G_{\alpha}'(u)\Delta u+G_{\alpha}''(u)|\nabla u|^2+G_{\alpha}'(u)\langle b,\nabla u\rangle-c\,G_{\alpha}(u)\\
&= G_{\alpha}''(u)|\nabla u|^2+G_{\alpha}'(u)\,c\,u-c\,G_{\alpha}(u)\,,
\end{aligned}
\end{equation}
where we have used \eqref{Laplacian_comp}, \eqref{soluzione_composta} and \eqref{problema}. Similarly to the proof of Proposition \ref{prop3}, by combining together \eqref{eq52} and \eqref{eq53b}, due to \eqref{G}, we get
\begin{equation}\label{eq53}
\begin{aligned}
\int_M &\left[G_{\alpha}'(u)\,c\,u-c\,G_{\alpha}(u)\right] G_{\alpha}(u)\gamma_R^2\,\phi^2\,d\mu\\
& \le -\left(1-\varepsilon \right) \int_M |G'_{\alpha}(u)|^2|\nabla u|^2\gamma_R^2\,\phi^2\,d\mu -  \int_M G''_{\alpha}(u)|\nabla u|^2\, G_{\alpha}(u)\,\gamma_R^2\,\phi^2\,d\mu \\
&\quad +\frac2{\varepsilon} \int_M  |\nabla \gamma_R|^2 G^2_{\alpha}(u)\,\phi^2\,d\mu + \frac2{\varepsilon} \int_M |\nabla \phi|^2 G^2_{\alpha}(u)\gamma_R^2\,d\mu\\
&\quad - \int_M \left\langle b\,,\,\nabla \gamma_R \right\rangle G^2_{\alpha}(u)\gamma_R\,\phi^2\,d\mu - \int_M \left\langle b\,,\,\nabla \phi\right\rangle G^2_{\alpha}(u)\gamma_R^2\,\phi\, d\mu\\
&\quad-\frac 12 \int_M\operatorname{div} b \,G_{\alpha}^2(u)\gamma_R^2\,\phi^2\,d\mu - \int_M c\,G_{\alpha}^2(u)\gamma_R^2\,\phi^2\,d\mu.\\
\end{aligned}
\end{equation}
On the left-hand side of \eqref{eq53}, by substituting the definition of $G_{\alpha}$ given in \eqref{G}, we have
\begin{equation}\label{eq53f}
\begin{aligned}
\int_M &\left[G_{\alpha}'(u)\,c\,u-c\,G_{\alpha}(u)\right] \,G_{\alpha}(u)\gamma_R^2\phi^2\,d\mu\\
&=\int_M (u^2+\alpha)^{\frac p4-1}\,c\left[\left(\frac p2-1\right)u^2-\alpha\right] (u^2+\alpha)^{\frac p4}\gamma_R^2\,\phi^2\,d\mu\,.
\end{aligned}
\end{equation}
Similarly, on the right-hand side of \eqref{eq53}, by exploiting the definition of $G_{\alpha}$ given in \eqref{G}, arguing as in \eqref{eq43g}, we can choose $\varepsilon>0$ sufficiently small such that, for any $p>1$,
\begin{equation}\label{eq53g}
-\int_M \left\{(1-\varepsilon)|G'_{\alpha}(u)|^2+G''_{\alpha}(u)G_{\alpha}(u)\right\}|\nabla u|^2\gamma_R^2\phi^2\,d\mu\le 0.
\end{equation}
Hence, \eqref{eq53}, combined with \eqref{eq53f} and \eqref{eq53g}, reduces to
\begin{equation}\label{eq53c}
\begin{aligned}
\int_M (u^2+\alpha)^{\frac p4-1}\,&c\left[\left(\frac p2-1\right)u^2-\alpha\right] (u^2+\alpha)^{\frac p4}\gamma_R^2\,\phi^2\,d\mu\\
&\le +\frac2{\varepsilon} \int_M  |\nabla \gamma_R|^2 G^2_{\alpha}(u)\,\phi^2\,d\mu + \frac2{\varepsilon} \int_M |\nabla \phi|^2 G^2_{\alpha}(u)\gamma_R^2\,d\mu\\
&\quad - \int_M \left\langle b\,,\,\nabla \gamma_R \right\rangle G^2_{\alpha}(u)\gamma_R\,\phi^2\,d\mu - \int_M \left\langle b\,,\,\nabla \phi\right\rangle G^2_{\alpha}(u)\gamma_R^2\,\phi\, d\mu\\
&\quad-\frac 12 \int_M\operatorname{div} b \,G_{\alpha}^2(u)\gamma_R^2\,\phi^2\,d\mu - \int_M c\,G_{\alpha}^2(u)\gamma_R^2\,\phi^2\,d\mu.\\
\end{aligned}
\end{equation}
Letting $\a\to 0^+$ in \eqref{eq53c}, by the dominated convergence and by choosing $\hat\delta=\frac 2\varepsilon$, we obtain the analogue of \eqref{eq44}, i.e.
\begin{equation}\label{eq54}
\begin{aligned}
\int_M |u|^p\,\gamma_R^2&\left[-\hat\delta |\nabla \phi|^2+\left\langle b\,,\,\nabla \phi\right\rangle\,\phi+\frac 12 ( \di b+ p\,c\,)\phi^2\right]\, d\mu\\
&\le \int_M |u|^p\left[\hat\delta |\nabla \gamma_R|^2\,\phi^2- \left\langle b\,,\,\nabla \gamma_R \right\rangle\gamma_R\,\phi^2\right]\,d\mu.
\end{aligned}
\end{equation}
Due to \eqref{eq54}, and the definitions of $D_+$ and $D_-$ in \eqref{d}, we obtain
\begin{equation}\label{eq56}
\begin{aligned}
&\int_M |u|^p\,\gamma_R^2\left[-\hat\delta |\nabla \phi|^2+\left\langle b\,,\,\nabla \phi\right\rangle\,\phi+\frac 12  (\di b\,+ p\,c)\,\phi^2\right]\, d\mu\\
&\quad \le\hat\delta \int_{M} |u|^p\,\phi^2 |\nabla \gamma_R|^2\,d\mu-\int_{D_+}|u|^p\,\phi^2 \left\langle b,\nabla \gamma_R \right\rangle\gamma_R\,d\mu -\int_{D_-}|u|^p \,\phi^2 \left\langle b,\nabla \gamma_R \right\rangle\gamma_R\,d\mu.
\end{aligned}
\end{equation}
Observing that, for all $x\in D_-$
$$
- \langle b(x)\,,\,\nabla \gamma_R(x)\rangle \,\le 0,
$$
we rewrite \eqref{eq56} as follows
\begin{equation}\label{eq56a}
\begin{aligned}
\int_M |u|^p\,\gamma_R^2&\left[-\hat\delta |\nabla \phi|^2+\left\langle b\,,\,\nabla \phi\right\rangle\,\phi+\frac 12  (\di b\,+ p\,c)\,\phi^2\right]\, d\mu\\
&\le \hat\delta\,\int_{M} |u|^p\,\phi^2 |\nabla \gamma_R|^2\,d\mu-\int_{D_+}|u|^p\,\phi^2 \left\langle b,\nabla \gamma_R \right\rangle\gamma_R\,d\mu.
\end{aligned}
\end{equation}
We now set $\phi(x):=e^{\zeta(x)}$ for any $x\in M$. Observe that, $\phi$ satisfies assumptions \eqref{eq52d} and \eqref{eq52c}. Hence, we can apply Lemma \ref{lemma5} with this choice of $\phi$. From Lemma \ref{lemma5} and the monotone convergence theorem, sending $R\to \infty$ in \eqref{eq56a} we get
\begin{equation}\label{eq57}
\int_M |u|^p \left[-\hat\delta\, |\nabla \phi|^2+\left\langle b\,,\,\nabla \phi\right\rangle\,\phi+\frac 12  (\di b+ p\,c)\,\phi^2\right]\, d\mu\leq 0\,.
\end{equation}
From \eqref{eq57} and \eqref{eq50}, since $|u|^p\geq 0$, we can infer that $u\equiv 0$ in $M$. This completes the proof.
\end{proof}

\begin{proof}[Proof of Theorem \ref{teo2}]
Let $\phi=\phi(r):=\xi(r)^{\frac 12}$ for any $r>0$, where $\xi$ has been defined in \eqref{eq13}, with $\tau>0$. 
At first observe that $\phi^2$ satisfies \eqref{eq52d} and \eqref{eq52c}. Moreover, $\phi$ solves \eqref{eq50}. Indeed, due to \eqref{h1} and \eqref{h2}, we have
\begin{equation}\label{eq55}
\begin{aligned}
&\hat\delta\, |\nabla \phi|^2-\left\langle b\,,\,\nabla \phi\right\rangle\,\phi-\frac 12  (\di b+ c\,p)\,\phi^2\\
&=\frac{\tau\hat\delta}{2}\left(\frac{\tau}{2}+1\right) (1+r)^{-\tau-4}|\nabla r|^2+\frac {\tau}2 (1+r)^{-\tau-1} \left\langle b\,,\,\nabla r\right\rangle -\frac 12 ( \di b+ c\, p)\,(1+r)^{-\tau}\\
&\le \frac{\tau\hat\delta}{4}(\tau+2) (1+r)^{-\tau-4}+\frac {\tau}2 \,K_2(1+r)^{-\tau-1+\s}+\frac {K_2}2(1+r)^{-\tau+\s-1}- \frac {c_0\, p}2 (1+r)^{-\tau}\\
&\le -\frac 12 (1+r)^{-\tau}\left(c_0\,p-\frac{\tau\hat\delta}{2}(\tau+2) (1+r)^{-4}- K_2\,(\tau+1)(1+r)^{\s-1}\right)\, . 
\end{aligned}
\end{equation}
Then, by \eqref{eq55}, one has
$$
\hat\delta |\nabla \phi|^2-\left\langle b\,,\,\nabla \phi\right\rangle\,\phi-\frac 12  (\di b+ c\,p)\,\phi^2<0,
$$ 
provided that
$$
p\,c_0>\frac{\tau\hat\delta}{2}(\tau+2)+K_2(\tau+1).
$$
Thus, by Proposition \ref{prop5} the conclusion follows. 
\end{proof}

\

\section{Proof of Proposition \ref{prop4} and Corollary \ref{cor}}\setcounter{equation}{0}\label{p3}

\

We now show nonuniqueness for problem \eqref{problema} with $b$ as in \eqref{h3} in the special case of $M$ being a model manifold satisfying \eqref{metric}. From standard results, see e.g. \cite[Theorem 2.5 and Proposition 2.7]{PTes}, the following Proposition can be proved.

\begin{proposition}\label{prop42}
Let $\mathbb{M}^N_\varphi$ be a model manifold such that \eqref{metric} and \eqref{vol2} are satisfied. Let $b:\mathbb{M}^N_\varphi\to \mathbb{M}^N_\varphi$, $b\in C^1(\mathbb{M}^N_\varphi)$ be such that assumption \eqref{h3} holds. Moreover assume \eqref{h2}. Suppose that there exist 
\begin{itemize}
\item[(i)] a bounded supersolution $h$ of equation 
\begin{equation}\label{eq100}
\Delta h+\left\langle b,\,\nabla h\right\rangle - c h =-1 \quad\quad \text{in}\,\,\, \mathbb{M}^N_\varphi\setminus B_{R_0}\,,
\end{equation}
for some $R_0>0$, such that 
\[h>0 \quad \text{ in }\,\, M, \quad \lim_{r\to +\infty} h(x)=0\,.\]
\item[(ii)] a positive bounded supersolution $W$ of the equation
\begin{equation}\label{eq101}
\Delta W+\left\langle b,\,\nabla W\right\rangle - c W =-1 \quad\quad \text{in}\,\,\, \mathbb{M}^N_\varphi\,.
\end{equation}
\end{itemize}
Then there exist infinitely many bounded solutions $u$ of problem \eqref{problema}. 
In particular, for any $\gamma\in \mathbb R$, there exists a solution $u$ to equation \eqref{problema} such that
$$
\lim_{r\to\infty} u(x)=\gamma.
$$
\end{proposition}

Proposition \ref{prop42} shows that infinitely many solutions to equation \eqref{problema} exist, if a supersolution $h$ to \eqref{eq100} exists and if a superolution $W$ to \eqref{eq101} exists. Therefore the following existence results, combined with the above proposition, imply nonuniqueness.

\begin{lemma}\label{lemma6bis} 
Let assumptions \eqref{metric}, \eqref{vol2}, \eqref{h2} and \eqref{h3} hold. Then there exists a supersolution $h>0$ of equation \eqref{eq100}.
\end{lemma}

\begin{proof} 
We consider the function
\begin{equation*}\label{eq540}
h(x):= Cr^{-\beta}\quad \text{for any}\,\,x\in \mathbb{M}^N_\varphi\setminus B_{R_0},
\end{equation*}
where $C>0$ and $\beta>0$ have to be chosen. 
Note that, due to \eqref{Laplaciano_modelli} and \eqref{metric}
$$
\begin{aligned}
\Delta h(x)&=\beta(\beta+1)C\,r^{-\beta-2} - \frac{\varphi'}{\varphi}(N-1)\,\beta\,C \,r^{-\b-1}\\
&=\beta C\left(\beta+1-(N-1)\,\lambda\right )\,r^{-\beta-2}\quad \quad \quad \text{for all}\,\,x\in \mathbb{M}^N_\varphi\setminus B_{R_0}.
\end{aligned}
$$
Thus, in view of \eqref{h3}, we have, for some $C_{\beta}>0$.
\begin{equation}\label{eq542}
\begin{aligned}
\Delta h(x)&+\langle b(x),\,\nabla h(x) \rangle-c\,h(x)\,\\
&\le \beta C\left(\beta+1-(N-1)\,\lambda\right )\, r^{-\beta-2}-\beta C r^{-\beta-1}\left\langle b(x),\,\nabla r\right\rangle-c C\,r^{-\b}\\
&\le \beta C\left(\beta+1-(N-1)\,\lambda\right )\, r^{-\beta-2}-\beta C K r^{\sigma-\beta-1}-cC\,r^{-\b} \\
&\le -C_{\b} \, r^{\sigma-1-\beta},
\end{aligned}
\end{equation}
for all $x\in \mathbb{M}^N_\varphi\setminus B_{R_0}$, with $R_0$ as in \eqref{h3}. Now, from \eqref{eq542} it follows that 
\[\Delta h(x)+\langle b(x),\,\nabla h(x) \rangle-c\,h(x)\,\le\,-1 \quad \text{for all}\,\,\,x\in \mathbb{M}^N_\varphi\setminus\bar{B}_{R_0};\]
provided that
$$
\sigma-1-\beta>0. 
$$
Hence, by the assumption $\sigma>1$, for $\beta$ small enough, we get the thesis.
\end{proof}

\begin{lemma}\label{lemma6}
Let assumption \eqref{h2} holds. Then there exists a bounded supersolution $W>0$ of equation \eqref{eq101}.
\end{lemma}

\begin{proof}
We define the function
$$
W(x):=\frac 1{c_0},
$$
for $c_0$ as in \eqref{h2}. Then, due to \eqref{h2}
\begin{equation*}\label{eq549}
\Delta W(x)+\langle b(x),\,\nabla W(x) \rangle-c\,W(x)\,\le -\frac{c}{c_0}\le -1,
\end{equation*}
for all $x\in \mathbb{M}^N_\varphi$. The Lemma is proved.
\end{proof}

\begin{proof}[Proof of Proposition \ref{prop4}]
By simply combining Proposition \ref{prop42}, Lemma \ref{lemma6bis} and Lemma \ref{lemma6}, the result follows.
\end{proof}

\bigskip

\begin{proof}[Proof of Corollary \ref{cor}]
\noindent \textit{(i)} Due to \eqref{metric}, it follows that, for any $B_R\subset M$,
$$
\operatorname{Vol}(x_0,r)\le C\,r^{(N-1)\lambda+1}.
$$
Consequently, $M$ satisfies assumption \eqref{vol2} for $\alpha=(N-1)\lambda+1$. Therefore, since $b$ satisfies \eqref{h1}, the thesis follows by Theorem \ref{teo2}.
\medskip

\noindent \textit{(ii)} The thesis directly follows by applying Proposition \ref{prop4}.
\end{proof}
\medskip

\

\noindent{\bf Acknowledgement} 

\noindent The authors have been partially supported by the ``Gruppo Nazionale per l'Analisi Matematica, la Probabilit\'a e le loro Applicazioni'' (GNAMPA) of the ``Istituto Nazionale di Alta Matematica'' (INdAM, Italy). A.R. has been partially supported by the GNAMPA project 2022 ``Stime ottimali per alcuni funzionali di forma''.

\

\

\noindent{\bf Data availability statement}

\noindent Data sharing not applicable to this article as no datasets were generated or analysed during the current study.

\

\end{document}